\documentclass[12pt,fullpage]{article}
\usepackage{amssymb,amsxtra,amsmath,amsthm}
\def \NN{{\mathbb{N}}}
\def \ZZ{{\mathbb{Z}}}
\def \QQ{{\mathbb{Q}}}
\def \RR{{\mathbb{R}}}

\def \FF{{\mathbb{F}}}
\def \PP{{\mathbb{P}}}
\def \CCC{{\mathcal{C}}}
\def \OOO{{\mathcal{O}}}
\def \TTT{{\mathcal{T}}}

\def \mm{{\mathfrak m}}
\def \cc{{\mathfrak c}}

\def\mod{\mathop{\mathrm{mod}}}

\def\det{\mathop{\mathrm{det}}}
\def\vert{\mathop{\mathrm{vert}}}
\def\edge{\mathop{\mathrm{edge}}}
\def\star{\mathop{\mbox{\Large$*$}}}

\def\Gamma{\varGamma}

\begin{document}
\begin{center} {\large {\bf NONRATIONAL GENUS ZERO FUNCTION FIELDS AND THE BRUHAT-TITS TREE}}\\
\bigskip
{\tiny {\rm BY}}\\
\bigskip
{\sc A. W. Mason}\\
\bigskip
{\small {\it Department of Mathematics, University of Glasgow\\
Glasgow G12 8QW, Scotland, U.K.\\
e-mail:awm@maths.gla.ac.uk}}\\
\bigskip
{\tiny {\rm AND}}\\
\bigskip
{\sc Andreas Schweizer}\\
\bigskip
{\small {\it National Center for Theoretical Sciences, Mathematics Division,\\
National Tsing Hua University\\
Hsinchu 300, Taiwan\\
e-mail: schweizer@math.cts.nthu.edu.tw}}
\end{center}
\begin{abstract}
Let $K$ be a function field with constant field $k$ and let $\infty$
be a fixed place of $K$. Let $\CCC$ be the Dedekind domain
consisting of all those elements of $K$ which are integral outside
$\infty$. The group $G=GL_2(\CCC)$ is important for a number of
reasons. For example, when $k$ is finite, it plays a central role in
the theory of Drinfeld modular curves. Many properties follow from
the action of $G$ on its associated {\it Bruhat-Tits tree}, $\TTT$.
Classical Bass-Serre theory shows how a presentation for $G$ can be
derived from the structure of the quotient graph (or fundamental
domain) $G\backslash \TTT$. The {\it shape} of this quotient graph
(for any $G$) is described in a fundamental result of Serre. However
there are very few known examples for which a detailed description
of $G\backslash \TTT$ is known. (One such is the {\it rational}
case, $\CCC=k[t]$, i.e. when $K$ has genus zero and $\infty$ has
degree one.) In this paper we give a precise description of
$G\backslash \TTT$ for the case where the genus of $K$ is zero, $K$
has no places of degree one and $\infty$ has degree two. Among the
known examples a new feature here is the appearance of vertex
stabilizer subgroups (of $G$) which are of {\it quaternionic} type.\\

\vspace{0.30cm} \noindent {\bf Mathematics Subject Classification
(2000):} 11F06, 20E08, 20H10, 11R58
\end{abstract}
\begin{center}
\section*{Introduction}
\end{center}
\noindent Let $K$ be an algebraic function field of one variable
over its full constant field $k$ and let $g$ be the genus of $K$.
From now on we fix a place $\infty$ of $K$. Let $d$ be the
degree of $\infty$. The principal focus of our attention will be
$\CCC$, the subset of $K$ consisting of all those elements of $K$
which are integral outside $\infty$. It is well-known that $\CCC$ is
a Dedekind domain and that $\CCC^*=k^*$. The group $G=GL_2(\CCC)$ is
important for a number of reasons. For example, when $k$ is finite,
it plays a central role [G] in the theory of Drinfeld modular
curves. The first intensive study of these groups can be
found in Serre's book [Se2].\\
\noindent Associated with $G$ is its {\it Bruhat-Tits tree}, $\TTT$.
(See [Se2, Chapter II, Section 1].) Classical Bass-Serre theory
applies to the action of $G$ on $\TTT$ and shows [Se2, Theorem 13,
p.55] how a presentation for $G$ can be obtained from the structure
of the quotient graph (or fundamental domain) $G \backslash \TTT$.
Of fundamental importance [Se2, Theorem 9, p.106] is the following
result.
\\ \\
\noindent {\bf Serre's Theorem.} {\it There exists a set $\Lambda$
in one-one correspondence with $Cl(\CCC)$, the ideal class group of
 $\CCC$, and a collection of pairwise-disjoint infinite paths, each
without backtracking, $\{X_{\lambda}:\lambda \in \Lambda\}$,
together with a connected graph $Y$ of bounded width, such that}
$$ G\backslash \TTT=
\begin {displaystyle} \left( \bigcup_{\lambda \in \Lambda}X_{\lambda}\right)
\cup Y,\end{displaystyle}$$
{\it where, for each $\lambda \in\lambda$,}
\begin{itemize}
\item[(i)] $|\mathrm{vert}(Y) \cap \mathrm{vert}(X_{\lambda})|=1$,
\item[(ii)] $|\mathrm{edge}(Y)\cap \mathrm{edge}(X_{\lambda})|=0.$
\end{itemize}
\noindent For the case where $k$ is finite Lubotzky [L, Theorem 6.1]
has extended this result to include all arithmetic lattices in rank
$1$ algebraic groups over nonarchimedean local fields. Serre's proof
is based on the theory of rank $2$ vector bundles over the
nonsingular projective curve corresponding to $K$. For a more
elementary alternative approach, which refers explicitly to
matrices, see [M1]. Although the stabilizers of all vertices and
edges of each $X_{\lambda}$ are well-understood not nearly as much
is known in general about $Y$. It is however known [Se2, Remark 3,
p.97] that one of the vertices of $Y$ has stabilizer $GL_2(k)$. It
is also known [Se2, Corollary 4, p.108] that $Y$ is finite if $k$ is
finite. In the case where $k$ is finite the theory of Drinfeld
modular curves provides an arithmetic formula [G, Chapter V] for the
rank of $\pi_1(Y)$, the fundamental group of (the graph) $Y$. (This
result has been used [MS], for example, to determine noncongruence
subgroups of $G$ of {\it minimal} index.) However the precise
structure of $Y$ is known in very few instances. To date they
consist of (i) the {\it rational} case, i.e. when $g=0$ and $K$ has
a place of degree $1$, [Se2, p.87], [M2], and (ii) the {\it
elliptic} case [T], i.e. when $g=d=1$. From the Riemann-Roch Theorem
it immediately follows that, if $K$ has genus zero, then the
smallest degree of any of its places is $1$ or $2$. The former
possibility gives rise to the rational case, i.e. $K=k(t)$. In this
paper we are concerned with the latter, {\it nonrational} case. By
virtue of a well-known result of F.K. Schmidt [St, V.1.12 Corollary,
p.164] the field of constants $k$ here is {\it infinite} (in which
case both $\CCC$ and $G$ are ``nonarithmetic"). Moreover it turns
out that nonrational function fields of this type
are subject to a number of restrictions.
\\
\noindent In our main result (Theorem 2.9) we determine completely
the structure of $G\backslash \TTT$, when $g=0$, $K$ is nonrational
and $d=2$. Of especial interest here is the appearance of vertex
stabilizer subgroups of $G$ which are {\it quaternion} groups. None
of the stabilizers of the vertices in any of the ``ends"
$X_{\lambda}$ of $G\backslash \TTT$ in Serre's theorem or the
Bruhat-Tits trees for the rational and elliptic cases are of this
type. We also determine, for several specific cases, the structure
of $\Gamma\backslash\TTT$, where $\Gamma=SL_2(\CCC)$. One important
feature which (sometimes) distinguishes $\Gamma$ from $G$ is the
existence of a free quotient. Among the more interesting
consequences of our results is the following.
\\ \\
\noindent {\bf Corollary.} {\it Let $\mathcal{Q}=\QQ[x, y]$, where
$x^2+y^2+1=0$. Then there exists an epimorphism $$ SL_2(\mathcal{Q})
\twoheadrightarrow F_{\omega}\;,$$ where $F_{\omega}$ is the free
group of countably infinite rank.
\par
In particular, every finitely or countably generated group
(including every finite group) is a homomorphic image of
$SL_2({\mathcal Q})$.}
\\ \\
\noindent A result like this, of course, depends very much on the
field of constants. If, for example, we replace $\QQ$ with $\RR$ it
collapses completely.  Less detailed versions of some of our results
are stated {\it without proof}  by Serre [Se2, p.114, 115, 118]. In
contrast to [Se2] however our approach here, as in [M1], [M2] and
[T], is elementary. Most of our results are explicitly stated in
terms of matrices. Results similar to our Corollary are already
known. For example Krstic and McCool [KM] have proved that if $R$ is
any domain which is not a field then $SL_2(R[x])$ has a free
quotient of infinite rank. (This extends previous results of
Grunewald, Mennicke and Vaserstein [GMV].)

\begin{center}
\subsection*{1. Nonrational genus zero function fields}
\end{center}

\noindent For convenience we list at this point the notation which
will be used in the rest of the paper:

\begin{tabular}{ll}
$K$         & a genus zero algebraic function field of nonrational type;\\
$k$& the constant field of $K$;\\
$\infty$    & a fixed place of $K$ of degree $2$;\\
$\nu_{\infty}$& the (discrete) valuation of $K$ given by $\infty$;\\
$\pi$& a local parameter at $\infty$ in $K$;\\
$\OOO$    & the valuation ring of $\infty$;\\
$\mm=(\pi)$ & the maximal ideal of $\OOO$;\\
$\CCC$      & the ring of all elements of $K$ that are integral
outside $\infty$;\\
$K_{\infty}$& the completion of $K$ with respect to $\infty$;\\
$\OOO_{\infty}$& the completion of $\OOO$ with respect to
$\infty$;\\
$\TTT$ & the Bruhat-Tits tree of $GL_2(K_\infty)$;\\
$G$ & the group $GL_2(\CCC)$;\\
$\Gamma$ & the group $SL_2(\CCC)$;\\
$S_w$ & the stabilizer in $G$ of $w \in \vert(\TTT)\cup\edge(\TTT)$.
\end{tabular}

\noindent Let $\PP_K$ denote the set of places of $K$. Then, by
definition,

$$ \CCC=\{z \in K:\nu_P(z) \geq 0,\;\forall P \in \PP_K,\;P \neq \infty \}.$$

\noindent Then $\nu_{\infty}(z) \leq 0$, for all $z \in \CCC$. For
each $n \in \NN\cup\{0\}$ we define
$$\CCC(n):=\{z \in \CCC: \nu_{\infty}(z) \geq -n\}.$$
\noindent In particular we have $\CCC(0)=k$. Although the characterization
of nonrational genus zero function fields is well-known
([A, Chapter 16, Section 4] and [VS, Section 9.6.1]),
we provide an outline of the proof, because a precise and explicit
description of $\CCC$ is essential for our purposes.
For simplicity we denote from now on $\nu_{\infty}$ by $\nu$.
\\ \\
\noindent {\bf Theorem 1.1.} \it Let $K/k$ be a nonrational function
field of genus zero. Then there exist $x,y \in K$, with
$\nu(x)=\nu(y)=-1$, such that
$$\CCC=k[x,y].$$
\noindent Furthermore $x,y$ can be chosen to satisfy equations of the
following types (with $\rho,\sigma \in k^*$).\\ \noindent
When $\mathrm{char}(k) \neq 2$
$$(i) \;\;\;y^2+\rho x^2+\sigma=0.$$
\noindent When $\mathrm{char}(k)=2$ there are three possibilities.
$$(ii)\;\;\;y^2+\rho x^2+\sigma=0.$$
$$(iii)\;\;\;y^2+\rho x^2+x+\sigma=0.$$
$$(iv)\;\;\;y^2+xy+\rho x^2+\sigma=0.$$
\noindent Moreover, the above equations have no solutions $x,y$ in $k$.
Finally, in cases $(i)$, $(ii)$ and $(iii)$ the element $-\rho$ must be
a non-square in $k$, whereas in case $(iv)$ $\rho$ is not of the form
$\rho=\alpha^2 +\alpha$ for any $\alpha\in k$.
\par
Conversely, if $\rho$ and $\sigma$ are elements of $k^*$ such that
any of the above equations has no nontrivial solutions $x,y$ in $k$,
and $-\rho$ is a non-square in $k^*$, resp. $\rho$ is not of the form
$\alpha^2 +\alpha$ in case $(iv)$,
then $k(x,y)$ is a non-rational function field of genus zero. \rm
\\ \\
\noindent {\bf Proof.} By the Riemann-Roch theorem [St, I.5.17, p.29]
$$\mathrm{dim}_k(\CCC(n))=2n+1.$$
In particular
$\mathrm{dim}_k(\CCC(1))=3$ and so there exist $x,y \in \CCC$ such
that $\nu(x)=\nu(y)=-1$ and $x,y,1$ are $k$-linearly independent.
Now $\mathrm{dim}_k(\CCC(2))=5$ and so
$x^2,y^2,xy,x,y,1$ are $k$-linearly dependent. Hence
$$\alpha y^2+\beta xy+ \gamma x^2+\delta x+\epsilon y+\zeta=0,$$
for some $\alpha,\beta,\gamma,\delta,\epsilon,\zeta \in k$, not all zero.
Now we may assume that $\alpha=1$. Then a $k$-basis for $\CCC(n)$ is
$$ 1,\cdots,x^n,y,yx,\cdots,yx^{n-1}.$$
It follows that
$\CCC=k[x,y]$. Thus $K$, the quotient field of $\CCC$, is a
quadratic extension of the rational function field $k(x)$.
\\
If $\mathrm{char}(k) \neq 2$ by completing squares in the above
equation we may assume that $\beta=\delta=\epsilon=0$. Now suppose
that $\mathrm{char}(k) = 2$. There are two possibilities. If
$\beta=0$ then, replacing $x$ or $y$ with $k$-linear combinations of
$x,y$ we can assume that $\epsilon=0$ and that $\delta =0,1$.
Finally if $\beta \neq 0$ we can replace $x,y$ with similar
combinations to reduce the equation to the case where $\beta=1$ and
$\delta=\epsilon=0$.
\\
By [A, p.302ff], in any characteristic the function field of a conic is
nonrational of genus $0$ if and only if the conic has no rational point
(in projective space). The further restrictions in Theorem 1.1 are necessary
to ensure this, or in other words, they insure that $K$ has no places of
degree $1$. For example the condition on $-\rho$ guarantees that the
infinite place of $k(x)$ is inert in $K/k(x)$.
\hfill $\Box$
\\ \\
{\bf Remark 1.2.}
Dividing the equation $(iv)$ by $x^2$ or by $y^2$ and changing variables,
we see that each field $K$ of type $(iv)$ can also be represented by
an equation of type $(iii)$, and vice versa. But this does not hold for
the rings $\CCC$. In case $(iii)$ the residue field of the place $\infty$
is inseparable over $k$, whereas for case $(iv)$ it is separable.
\par
Since $K$ has no places of degree $1$, all places of degree $1$ in $k(y)$
must be inert in $K/k(y)$. If $K$ is separable over $k(y)$, as it is in
cases $(iii)$ and $(iv)$, no place of $k(y)$ can be ramified in $K/k(y)$
since the wildness of the ramification would force $K$ to have positive
genus. But by the Hurwitz formula the different of $K/k(y)$ must have
degree $2$. Thus exactly one place of degree $1$ of $k(y)$ has an inseparable
residue extension in $K$. In case $(iii)$ this place is $\infty$, whereas
in case $(iv)$ it is the place lying above $y=0$.
\\ \\
As one might guess from Theorem 1.1, not every field $k$ allows the
existence of a nonrational genus zero function field $K/k$. We list
the most important instances.
\\ \\
{\bf Remark 1.3.} Let $F$ be a function field of genus zero with
constant field $k$. In the following situations, $F$ is automatically
a rational function field:
\begin{itemize}
\item[a)] $k$ is algebraically closed (clear);
\item[b)] $k$ is a finite field (This follows from the Theorem of
F.K.~Schmidt. Alternatively: the pigeon-hole principle shows that over
a finite field the equations from Theorem 1.1 always have a solution.);
\item[c)] $k$ is any algebraic extension of a finite field (Using
Theorem 1.1 this is easily reduced to the case where $k$ is finite);
\item[d)] $k$ itself is a function field of one variable over an
algebraically closed constant field (In this case $k$ is a $C_1$-field.
This means that every homogenous polynomial over $k$ of degree $d$ in
more than $d$ variables has a non-trivial zero. We only need the case
$d=2$. See [FJ, Proposition 21.2.12 for a more general statement).
\end{itemize}
In characteristic $2$ one can completely characterize the fields $k$ over
which case $(ii)$ of Theorem 1.1 can occur.
We write $k^2$ for the field consisting of all squares in $k$. Note that
$k^2$ is isomorphic to $k$ via the Frobenius. The degree $[k:k^2]$ is
a power of $2$ (possibly infinite). For example, if $L$ is perfect and
$k=L(x_1,x_2,\ldots,x_n)$ where the $x_i$ are algebraically independent,
then $k^2 =L(x_1^2,x_2^2,\ldots,x_n^2)$ and $[k:k^2]=2^n$.
\\ \\
{\bf Lemma 1.4.} \it
In characteristic $2$, case $(ii)$ of Theorem 1.1 occurs if and only if
$[k:k^2]>2$; i.e. if and only if $k$ has at least two inseparable extensions
of degree $2$.
In particular, case $(ii)$ cannot occur in the following situations:
\begin{itemize}
\item[a)] $k$ is perfect;
\item[b)] $k$ itself is a function field of one variable
over a perfect constant field $L$;
\item[c)] $k$ is a Laurent series field over a perfect field $L$.
\end{itemize}
\rm
\noindent {\bf Proof.}
From Theorem 1.1 we know that $\rho$ must be a non-square. If
$[k:k^2]\le 2$, then $k^2 +\rho k^2=k$, so $y^2 +\rho x^2$ represents
$\sigma$ and the equation $(ii)$ from Theorem 1.1 has a solution.
If $[k:k^2]>2$, we can choose non-squares $\rho$ and $\sigma$ such
that $\sigma\not\in k^2 +\rho k^2$, so the equation has no solutions.
\par
a) is already clear from Theorem 1.1.
\par
b) Since $L$ is perfect, we can find an element $u$ in $k$ such that $k$
is a finite separable extension of $L(u)$. Hence there exists a primitive
element $v$ with $k=L(u,v)$. Since $L$ is perfect and sums of squares are
squares in characteristic $2$, the extension obtained by extracting a
square-root from every element of $k$ is
$$L(\sqrt{u},\sqrt{v}).$$
Let $f(X)\in L(u)[X]$ be the minimal
polynomial of $v$ over $L(u)$. Then $\sqrt{v}$ is a zero of $f(X^2)$,
which is the square of a polynomial in $L(\sqrt{u})[X]$. Thus
$$[L(\sqrt{u},\sqrt{v}):L(u)]
=[L(\sqrt{u}):L(u)]\cdot[L(\sqrt{u},\sqrt{v}):L(\sqrt{u})]=2[k:L(u)].$$
So $L(\sqrt{u},\sqrt{v})$ has degree $2$ over $k$. Thus it is the
only inseparable extension of degree $2$.
\par
c) One easily checks that the squares in $L((u))$ form the field $L((u^2))$,
so $[k:k^2]=2$.
\hfill $\Box$
\\ \\
{\bf Example 1.5.}
Lemma 1.4 shows in particular that over a local or global constant field $k$
of characteristic $2$ the case $(ii)$ of Theorem 1.1 cannot occur. However,
case $(iii)$ and hence also case $(iv)$ occur. For example, if the
quadratic polynomial $z^2 +z+\kappa$ is irreducible over $\FF_{2^r}$,
then the equation
$$y^2 +\frac{1}{u}x^2 +x+\kappa u=0$$
satisfies the conditions of Theorem 1.1 over $k=\FF_{2^r}((u))$
and hence also over $k=\FF_{2^r}(u)$ and over any finite extension of
$\FF_{2^r}(u)$ contained in $\FF_{2^r}((u))$.
\\

\begin{center}
{\subsection*{2. The fundamental domain}}
\end{center}

\noindent We opt for the model of $\TTT$ used by Takahashi [T]. (The
model used by Serre [Se2, Chapter II, 1.1] is different but
equivalent.) By definition
$$\OOO=\{ z \in K:\nu(z) \geq 0\}.$$
\noindent The {\it vertices} of $\TTT$ are the left cosets of
$ZGL_2(\OOO_{\infty})$ in $GL_2(K_{\infty})$, where
$ZGL_2(\OOO_{\infty})$ is the subgroup generated by
$GL_2(\OOO_{\infty})$ and $Z$, the centre (i.e. the scalar matrices)
of $GL_2(K_{\infty})$. We recall that $\nu(\pi)=1$. The {\it edges}
of  $\TTT$ are defined in the following way. The vertices
$g_1ZGL_2(\OOO_{\infty})$ and $g_2ZGL_2(\OOO_{\infty})$ are {\it
adjacent} if and only if
$$g_2^{-1}g_1 \equiv  \left[\begin{array}{lll} \pi & z\\[10pt]
0 & 1
\end{array}\right]\; \mathrm{or}\: \left[\begin{array}{lll} \pi^{-1} & 0\\[10pt]
0 & 1
\end{array}\right]\; (\mod ZGL_2(\OOO_{\infty})).
$$
for some $z \in K$. Then $\TTT$ is a tree on which $G$ acts (by left
multiplication of cosets) without inversion. \noindent Since
$\OOO_{\infty}$ is a PID we may represent every coset of
$ZGL_2(\OOO_{\infty})$ by an element of the form
$$\left[\begin{array}{lll} \pi^n & z\\[10pt]
0 & 1
\end{array}\right],
$$
for some $ n \in \ZZ$ and $z \in K_{\infty}$. We denote this vertex
of $\TTT$ by $v(\pi^n,z)$. It is clear that
$$ v(\pi^n,z)=v(\pi^m,z') \Longleftrightarrow n=m \;\mathrm{and}\; \nu(z-z') \geq n.$$
We may assume therefore that $ z \in K$. Let $S(\pi^n,z)$ denote the
stabilizer of $v(\pi^n,z)$ in $\Gamma$. We record the following
well-known result.

\noindent {\bf Lemma 2.1.} \it With the above notation, let
$$ g=\left[\begin{array}{lll} a & b\\[10pt]
c & d
\end{array}\right].$$
Then $ g \in S(\pi^n,z)$ if and only if
\begin{itemize}
\item[(i)]$ \nu(a-zc) \geq 0,\;\nu(d+zc)\geq 0,$
\item[(ii)] $\nu(c) \geq -n,$
\item[(iii)] $\nu(b+z(a-d)-z^2c)\geq n.$
\end{itemize}

\noindent \rm If $ v \in \vert(\TTT)$, then, $S_v$, the stabilizer
of $v$ in $G$, acts on the edges of $\TTT$ incident with $v$. Let
$\mathrm{Orb}(v)$ denote a complete set of representatives for this
action. Our determination of $G\backslash \TTT$ depends on the
following easily verifiable result.

\noindent {\bf Lemma 2.2.} \it Let $\TTT_0$ be the following
infinite path in $\TTT$
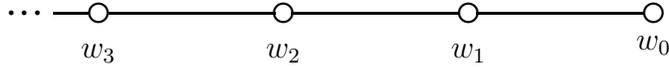
\begin{center}
\setlength{\unitlength}{1pt} \thicklines
\begin{picture}(300,50)(-150,-30)
\put(-110,0){\circle{7}} \put(-40,0){\circle{7}}
 \put(30,0){\circle{7}} \put(100,0){\circle{7}}
\put(-106,0){\line(1,0){63}} \put(-37,0){\line(1,0){63}}
 \put(33,0){\line(1,0){63}}
\put(-127,0){\line(1,0){13}} \put(-133,0){\circle{.5}}
\put(-138,0){\circle{.5}} \put(-143,0){\circle{.5}}
\put(-110,-15){\makebox(0,0){$w_3$}}
\put(-40,-15){\makebox(0,0){$w_2$}}
\put(30,-15){\makebox(0,0){$w_1$}}
\put(100,-12){\makebox(0,0){$w_0$}}
\end{picture}
\end{center}
\noindent and let $e_i$ be the edge of $\TTT$ joining $w_i$ and
$w_{i+1}$, where $i \geq 0$. \\ \\ \noindent Suppose that
\begin{itemize}
\item[(i)] $\mathrm{Orb}(w_0)=\{e_0\}$ and
$\mathrm{Orb}(w_i)=\{e_i,e_{i-1}\}$, where $i \geq 1$,
\item[(ii)]for all $i \geq 0$
$$S_{w_{i}} \ncong S_{w_{i+1}}\; .$$
 \end{itemize} Then
$$ G\backslash \TTT \cong \TTT_0.$$

\noindent \rm The proof of this lemma is a standard argument which
uses the fact that, if $u,v$ are vertices of $\TTT$ and $v=g(u)$,
for some $g \in G$, then $S_u$
and $S_v$ are conjugate in $G$. \\ \\
\noindent {\bf Notation.} Let $t=y/x$, where $x,y$ are defined as in
Theorem 1.1. We put
$$ v_n =v(\pi^{-n},0),$$ where $n \geq 0$, and
$$v_*=v(\pi,t).$$
\noindent Let $\epsilon_*$ be the edge of $\TTT$ joining $v_*$ and
$v_0$ and let $\epsilon_i$ be the edge of $\TTT$ joining $v_i$ and
$v_{i+1}$, where $i \geq 0$.\\ \\
\noindent In the remainder of this section we will show that the
infinite path in $\TTT$ whose vertices are $v_*$ and $v_n\;(n \geq
0)$ is isomorphic to the
fundamental domain $G\backslash \TTT$.\\ \\
\noindent Our next result is an easy consequence of Lemma 2.1.

\noindent {\bf Lemma 2.3.} \it
\begin{itemize}
\item[(i)] $S_{v_0}=S(1,0) = GL_2(k).$
\item[(ii)] For any $n \geq 1$,
$$S_{v_n}=S(\pi^{-n},0)=\left\{ \left[\begin{array}{lll} \alpha & b\\[10pt]
0 & \beta \end{array}\right]: \alpha,\beta \in k^*,\;b \in
\CCC(n)\right\}.$$
\end{itemize}

\noindent {\bf Lemma 2.4.} \it \begin{itemize}
\item[(i)] $\mathrm{Orb}(v_0)=\{\epsilon_*,\epsilon_0\}.$
\item[(ii)] $\mathrm{Orb}(v_n)=\{\epsilon_n,\epsilon_{n-1}\}\;\;(n
\geq 1).$
\end{itemize}

\noindent {\bf Proof.} \rm The fact that $\infty$ has degree $2$
means that the residue class field $\overline{k}=\OOO/\mm$ is a
quadratic extension of $k$. In the notation of Theorem 1.1,
$\overline{k}\cong k[\sqrt{-\rho}]$ in cases $(i)$, $(ii)$ and $(iii)$,
and $\overline{k}\cong k[\overline{t}]$ with
$\overline{t}^2 +\overline{t}+\rho=0$ in case $(iv)$.
\\ \\
\noindent From the above the vertices of $\TTT$ which are adjacent
to $v(1,0)$ are precisely $v(\pi^{-1},0)$ and $v(\pi,u)$, where $u
\in \overline{k}$. Now the action of $S_{v_0}$ on these is
equivalent to the action of $GL_2(k)$ on
$\overline{k}\cup\{\infty\}$, as a group of linear fractional
transformations. (See [Se2, Exercise 6, p.99].) The latter action
has precisely $2$ orbits, one containing $``{\infty}"$,
corresponding to $v(\pi^{-1},0)$, and the other $t$,
corresponding to $v_*=v(\pi,t)$. Part (i) follows.
\\ \\
\noindent When $n \geq 1$ the vertices of $\TTT$ which are adjacent
to $v_n=v(\pi^{-n},0)$ are precisely $v_{n+1}=v(\pi^{-(n+1)},0)$ and
$v(\pi^{-n+1},u\pi^{-n})$, where $u \in \overline{k}$. Now by Lemma
2.3 (ii) one of the orbits of $S_n=S(\pi^{-n},0)$ is
$\{\epsilon_n\}$. We now show that all the other edges of $\TTT$
incident with $v_n$ are in
the same $S_n$-orbit.\\ \\
\noindent By the Riemann-Roch theorem ($K$ has genus zero) it
follows that
$$\dim_k(\CCC(n+1)/\CCC(n))=2,$$
and hence that
$$\pi^n\CCC(n)+\mm=\OOO.$$ Let $u \in \overline{k}$. From the above
we may choose $b \in \CCC(n)$ such that
$$b\pi^n+u \equiv 0\;(\mod \pi).$$ It is easily verified that
$$ \left[\begin{array}{lll} \pi^{n-1} & 0\\[10pt]
0 & 1 \end{array}\right] \left[\begin{array}{lll} 1 & b\\[10pt]
0 & 1 \end{array}\right]\left[\begin{array}{lll} \pi^{-n+1} & u\pi^{-n}\\[10pt]
0 & 1 \end{array}\right] \in GL_2(\OOO).$$ \noindent Hence we may
represent $\mathrm{Orb}(v_n)$ by a subset of
$\{\epsilon_{n-1},\epsilon_n\}$. If these edges are in the same
$S_n$-orbit, then $g(v_{n-1})=v_n$, for some $g \in \Gamma$. It
follows that $S_n$ and $S_{n-1}$ are conjugate and hence isomorphic,
which contradicts Lemma 2.2. Part (ii) follows.\hfill $\Box$

\noindent There remains the most difficult vertex to deal with. Our
description is based on the notation used in Theorem 1.1.
\\ \\
\noindent {\bf Lemma 2.5.} \it With the notation of Theorem 1.1
$(i)$, $(ii)$ and $(iii)$, suppose that
$$y^2+\rho x^2+\tau x+\sigma=0,$$
where either $\tau=0$ or ($\tau=1$ and $\mathrm{char}(k)=2$). Let
$$U=\left[\begin{array}{cc} 0 &-\rho\\[10pt]
1 & 0 \end{array}\right],\;V= \left[\begin{array}{cc} y &\tau+\rho x\\[10pt]
x & -y \end{array}\right],\; W=\left[\begin{array}{cc} -\rho x & \rho y\\[10pt]
y & \tau+\rho x \end{array}\right].$$ \noindent Then
$$S(\pi,t)=\left\{\alpha I_2+\beta U+\gamma V+\delta
W:\alpha,\beta,\gamma,\delta \in
k,\;(\alpha,\beta,\gamma,\delta)\neq(0,0,0,0)\right\}.$$
\rm 
\noindent {\bf Proof.} \rm We recall from Theorem 1.1 that here
$\CCC=k[x,y]$, where $\nu(x)=\nu(y)=-1$. In addition $\nu(t)=0$ and
$t^2=-\rho-\tau x^{-1}-\sigma x^{-2}.$ Let
$$X=\left[\begin{array}{cc} a &b\\[10pt]
c& d \end{array}\right] \in S(\pi,t).$$ \noindent Then from Lemma
2.1
$$ \nu(c) \geq -1,\;\nu(a-tc),\nu(d+tc)\geq 0,\;\mathrm{and}\;
\nu(b+(a-d)t-ct^2) \geq 1. $$ \noindent We deduce that
$\nu(a),\nu(b),\nu(d)\geq -1$ as well. From the proof of Theorem 1.1
the elements $1,x,y$ form a $k$-basis for $\CCC(1)$, in which case
$$a=\alpha_a+\beta_ax+\gamma_a y,$$
\noindent where $\alpha_a,\beta_a,\gamma_a \in k$. There are
corresponding expressions for $b,c,d$.\\ \\
\noindent Now $\nu(a-tc) \geq 0$ and so $\nu((a-tc)x^{-1})\geq 1$.
Using the fact that $x^{-1} \in \mm$, so that $t^2 \equiv
-\rho\;(\mod \mm)$, it follows that
$$ (\beta_a+\gamma_c \rho)+(\gamma_a-\beta_c)t \equiv 0 \;(\mod
\mm).$$ \noindent Since $1,t$ form a $k$-basis for $\OOO/\mm$, we
deduce that
$$\beta_a=-\gamma_c \rho\;\; \mathrm{and}\;\;\gamma_a= \beta_c.$$
\noindent In an identical way it can be shown that
$$\beta_d=\gamma_c \rho \;\; \mathrm{and} \;\; \gamma_d=-\beta_c.$$
\noindent From the remaining condition satisfied by the entries of
$X$ it follows that $\nu((b+(a-d)t-ct^2)x^{-1}) \geq 2$. Using the
fact that $t^2 \equiv -\rho-\tau x^{-1}\;(\mod \mm^2)$ we deduce
from the above that
$$(\beta_b-\beta_c\rho)+(\alpha_b-\beta_c\tau+\alpha_c\rho)x^{-1}+
(\gamma_b-\gamma_c\rho)t+(\alpha_a-\alpha_d+\gamma_c\tau)tx^{-1}
\equiv 0 \;(\mod \mm^2).$$ \noindent Since $1,t,x^{-1},tx^{-1}$ form
a $k$-basis for $\OOO/\mm^2$ we conclude that
$$\alpha_d=\alpha_a+\gamma_c\tau,\;\alpha_b=\beta_c\tau-\alpha_c \rho,\;
\beta_b=\beta_c \rho,\;\gamma_b=\gamma_c \rho.$$ \noindent Then $X$
is of the form
$$\alpha I_2+\beta U+\gamma V+\delta W.$$
It is easily verified that
$$\det(X)=\alpha^2+\rho\beta^2+\tau(\alpha\delta-\beta\gamma)+\sigma(\gamma^2+\rho\delta^2).$$
\noindent Suppose that $\det(X)=0$. If $\gamma^2+\rho\delta^2=0$
then by Theorem 1.1 $\alpha=\beta=\gamma=\delta=0$. Assume now that
$\mu=\gamma^2+\rho\delta^2\neq 0$. It is easily verified that $$
\mu^{-1}\det(X)=\lambda_0^2+\rho\mu_0^2+\tau\mu_0+\sigma=0,$$ where
$\lambda_0=(\alpha\gamma+\rho\beta\delta)\mu^{-1}$ and
$\mu_0=(\alpha\delta-\beta\gamma)\mu^{-1}$. This contradicts Theorem
1.1. The rest of the proof is trivial.
\hfill $\Box$

\noindent We now deal with case $(iv)$ in Theorem 1.1.
\\ \\
\noindent {\bf Lemma 2.6.} \it With the notation of Theorem 1.1
$(iv)$, suppose that $\mathrm{char}(k)=2$ and that $$y^2+xy+\rho
x^2+\sigma=0.$$  Let
$$U=\left[\begin{array}{cc} 0 &\rho\\[10pt]
1 & 1 \end{array}\right],\;V= \left[\begin{array}{cc} y &\rho x+y\\[10pt]
x & y \end{array}\right],\; W=\left[\begin{array}{cc} \rho x+y & \rho x+(1+\rho) y\\[10pt]
y & \rho x+y \end{array}\right].$$ \noindent Then
$$S(\pi,t)=\left\{\alpha I_2+\beta U+\gamma V+\delta
W:\alpha,\beta,\gamma,\delta \in
k,\;(\alpha,\beta,\gamma,\delta)\neq(0,0,0,0)\right\}.$$
\rm
\noindent {\bf Proof.} In this case $t^2+t+\rho \equiv 0\;(\mod
\mm^2).$ The proof which is very similar to that of Lemma 2.5
follows from Theorem 1.1 and Lemma 2.1.
\hfill $\Box$
\\ \\
{\bf Remark 2.7.} Consider the case dealt with in Lemma 2.5 with $\tau=0$. Let
$$D=\{\alpha I_2+\beta U+\gamma V+\delta W:
\alpha,\beta,\gamma,\delta \in k\}.$$
\noindent Note that if $\{U,V,W\}=\{P,Q,R\}$ then $P^2$ is a scalar matrix
and $PQ=-QP$. Thus, if the characteristic of $k$ is different from $2$,
then $S(\pi,t)=D^*$, the non-zero elements of a non-split quaternion
algebra $D$ over $k$. If the characteristic is $2$ (i.e. case $(ii)$
from Theorem 1.1) however, $S(\pi,t)$ is isomorphic to the multiplicative
group of the field $k(\sqrt{\rho},\sqrt{\sigma})$.
\par
Similarly, in the cases $(iii)$ and $(iv)$ from Theorem 1.1 one sees
from Lemma 2.5 resp. Lemma 2.6 that the stabilizer $S(\pi,t)$ is the
multiplicative group $D^*$ of a quaternion skew-field $D$ over $k$.

We refer to $S(\pi,t)$ as a group of {\it quaternionic} type. It is
of especial interest since no vertex stabilizer in any $G$ for which
(i) $K$ is a genus zero function field of rational type or (ii) $K$
is elliptic (i.e. of genus $1$ with a rational point) is of this
form. See [M2], [T]. Serre's Theorem also shows that, for any $G$,
if $v$ is a vertex of $\TTT$ more than a bounded distance from, say,
the "standard" vertex $S(1,0)$, then $S_v$ is a nilpotent group (of
class at most $2$) which is not of quaternionic type.
\\ \\
\noindent {\bf Lemma 2.8.} \it With the above notation,
$$\mathrm{Orb}(v_*)=\{\epsilon_*\}.$$

\noindent {\bf Proof.} \rm The vertices of $\TTT$ which are adjacent
to $v_*=v(\pi,t)$ are $v(1,0)$ and $v(\pi^2,t+\pi u)$, where $u \in
\OOO/\mm$. It suffices therefore to prove that, for all $u \in
\OOO/\mm$, there exists $X \in S(\pi,t)$ such that
$$ X \equiv \left[\begin{array}{cc} \pi^2 &t+\pi u\\[10pt]
0 & 1 \end{array}\right] \;(\mod ZGL_2(\OOO_{\infty})).$$ \noindent
For simplicity we take $\pi=x^{-1}$. There are essentially two cases
to consider.
\\ \\
\noindent {\bf Case A.} With the notation of Theorem 1.1 suppose
that $$y^2+\rho x^2+\tau x+\sigma=0,$$ where either $\tau=0$ or
($\tau=1$ and $\mathrm{char}(k)=2$). As in Lemma 2.5 let
$$ V= \left[\begin{array}{cc} y &\tau+\rho x\\[10pt]
x & -y \end{array}\right],\; W=\left[\begin{array}{cc} -\rho x & \rho y\\[10pt]
y & \tau+\rho x \end{array}\right].$$
\noindent It is easily verified that
$$x^{-1}V\left[\begin{array}{cc} t &1\\[10pt]
1& 0 \end{array}\right] =\left[\begin{array}{cc}-\sigma x^{-2} & t\\[10pt]
0 & 1 \end{array}\right]\equiv  \left[\begin{array}{cc} \pi^2 &t\\[10pt]
0& 1 \end{array}\right] \;(\mod ZGL_2(\OOO_{\infty}))$$ \noindent
Now let
$$X=I_2+\alpha V+\beta W,$$
where $\alpha,\beta \in k$, with $(\alpha,\beta) \neq (0,0)$. Let
$z=(\alpha x+\beta y)$ and let $s=-z(1+\alpha y-\beta \rho x)^{-1}$.
Then $\nu(s)=0$ and it is easily verified that
$$z^{-1}X\left[\begin{array}{cc} 1 & 1\\[10pt]
s & 0 \end{array}\right]= \left[\begin{array}{cc} a' & b'\\[10pt]
0 & 1 \end{array}\right],$$ \noindent where (by considering
determinants) $a'=v\pi^2$, for some $v \in \OOO^*$, and
$$b'-t=\pi (1+\beta\tau+\beta\sigma x^{-1})(\alpha+\beta t)^{-1},$$
so that $$z^{-1}X\left[\begin{array}{cc} 1 & 1\\[10pt]
s & 0 \end{array}\right]\equiv \left[\begin{array}{cc} \pi^2 & \pi(1+\beta\tau)(\alpha +\beta t)^{-1}\\[10pt]
0 & 1 \end{array}\right]\;(\mod ZGL_2(\OOO_{\infty})).$$ \noindent
The proof of Case A is therefore complete apart from the case where
$\tau=\beta=1$. To deal with this consider
$$ Y = \alpha V+W.$$ Let $z=(\alpha x+y)$ and let
$s'=-z(\alpha y-\rho x)^{-1}$. Then $\nu(s')=0$ and it is easily
verified that
$$z^{-1}X\left[\begin{array}{cc} 1 & 1\\[10pt]
s' & 0 \end{array}\right]\equiv \left[\begin{array}{cc} \pi^2 & \pi(\alpha +t)^{-1}\\[10pt]
0 & 1 \end{array}\right]\;(\mod ZGL_2(\OOO_{\infty})).$$
\\
\noindent {\bf Case B.} Again with the notation of Theorem 1.1
suppose that $\mathrm{char}(k)=2$ and that$$y^2+xy+\rho
x^2+\sigma=0.$$
\noindent As in Lemma 2.6 let $$V= \left[\begin{array}{cc} y &\rho x+y\\[10pt]
x & y \end{array}\right],\; W=\left[\begin{array}{cc} \rho x+y & \rho x+(1+\rho) y\\[10pt]
y & \rho x+y \end{array}\right].$$ \noindent It is easily verified
that $$x^{-1}V \left[\begin{array}{cc} t &1\\[10pt]
1& 0 \end{array}\right]\equiv  \left[\begin{array}{cc} \pi^2 &t\\[10pt]
0& 1 \end{array}\right] \;(\mod ZGL_2(\OOO_{\infty})).$$ \noindent
Now let $$ X=I_2+\alpha V+\beta W,$$ where $(\alpha,\beta)\neq
(0,0)$. As before we put $z=\alpha x+\beta y$. Then as in the above
proof for Case A it can be shown that there exists $s$, with
$\nu(s)=0$, such that $$z^{-1}X\left[\begin{array}{cc} 1 & 1\\[10pt]
s & 0 \end{array}\right]\equiv \left[\begin{array}{cc} \pi^2 & \pi(\alpha +\beta t)^{-1}\\[10pt]
0 & 1 \end{array}\right]\;(\mod ZGL_2(\OOO_{\infty})).$$ \noindent
This completes the proof.
\hfill $\Box$

\noindent Combining Lemmas 2.1, 2.2, 2.3, 2.4, 2.5, 2.6 and 2.8 we
obtain the principal result in this paper.
\\ \\
\noindent {\bf Theorem 2.9.} \it With the above notation, the
infinite path in $\TTT$

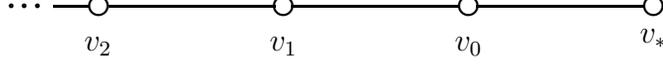
\begin{center}
\setlength{\unitlength}{1pt} \thicklines
\begin{picture}(300,50)(-150,-30)
\put(-110,0){\circle{7}} \put(-40,0){\circle{7}}
 \put(30,0){\circle{7}} \put(100,0){\circle{7}}
\put(-106,0){\line(1,0){63}} \put(-37,0){\line(1,0){63}}
 \put(33,0){\line(1,0){63}}
\put(-127,0){\line(1,0){13}} \put(-133,0){\circle{.5}}
\put(-138,0){\circle{.5}} \put(-143,0){\circle{.5}}
\put(-110,-15){\makebox(0,0){$v_2$}}
\put(-40,-15){\makebox(0,0){$v_1$}}
\put(30,-15){\makebox(0,0){$v_0$}}
\put(100,-12){\makebox(0,0){$v_*$}}
\end{picture}
\end{center}
is isomorphic to $G\backslash \TTT$.

\noindent \rm A group-theoretic consequence of Theorem 2.9 is the
following.
\\ \\
\noindent {\bf Corollary 2.10.} \it $G = A\star_{C}B,$ where
\begin{itemize}
\item[(i)]$A=S(\pi,t)$,
\item[(ii)] $B=GE_2(k[x,y])$, the subgroup of $G$ generated by
all its elementary matrices, \item[(iii)] $C=A\cap B \cong
(\OOO/\mm)^*.$ \end{itemize}

\noindent {\bf Proof.} \rm By the fundamental theorem of the theory
of groups acting on trees [Se2, Theorem 13, p.55]
 $$ G=A \star_{C} B,$$
 where \begin{itemize} \item[(i)] $A=S_{v_*}=S(\pi,t)$,
 \item[(ii)] $B=\langle S_{v_n}: n \geq 0\rangle=GE_2(k[x,y])$,
 \item[(iii)] $ C=A\cap B=\left\{\alpha I_2+\beta U: \alpha,\beta
 \in k^*,\; (\alpha,\beta) \neq (0,0)\right\} \cong (\OOO/\mm)^*$.
 \end{itemize} \hfill $\Box$
 \bigskip

\begin{center}
\subsection*{3. Subgroups with free quotients}
\end{center}

\noindent We recall that $\Gamma=SL_2(\CCC)$. We will identify
$G\backslash\TTT$ with the infinite path in $\TTT$ as in Theorem
2.9. The structure of $\Gamma\backslash \TTT$ can be derived from
that of $G\backslash \TTT$. In this way we will provide examples of
$\Gamma$ which have free quotients of up to
countably infinite rank.\\ \\
\noindent There is a natural epimorphism of graphs $$
\phi:\;\Gamma\backslash\TTT\twoheadrightarrow G\backslash\TTT.$$
\noindent For each $ w \in \vert(\Gamma\backslash
\TTT)\cup\edge(\Gamma\backslash\TTT)$ there is a one-one
correspondence
$$\phi^{-1}(w)\longleftrightarrow G/(S_w\cdot\Gamma),$$
the coset space of $S_w\cdot\Gamma$ in $G$.
\\ \\
Recall that the residue field of the place $\infty$ is $k(\overline{t})$
where $\overline{t}^2 +\rho=0$ in cases $(i)$, $(ii)$ and $(iii)$ of
Theorem 1.1 and $\overline{t}^2 +\overline{t}+\rho=0$ in case $(iv)$.
Hence from Lemmas 2.5 and 2.6 we obtain:
\par
If the characteristic of $k$ is different from $2$, the vertices of
$\Gamma\backslash\TTT$ lying above $v_*$ are in bijection with the
cosets of $N_{red}(D^*)$ in $k^*$ where $N_{red}$ is the reduced norm on the
quaternion algebra $D$. The edges above $\epsilon_*$ are in
bijection with the cosets
$k^*/N_{k(\sqrt{-\rho})/k}(k(\sqrt{-\rho})^*)$.
\par
If $\mathrm{char}(k)=2$,  the edges of $\Gamma\backslash\TTT$
above $\epsilon_*$ are in bijection with
the cosets $k^*/N_{k(\overline{t})/k}(k(\overline{t})^*)$. In cases $(iii)$
and $(iv)$ the vertices of $\Gamma\backslash\TTT$ lying above $v_*$ are in
bijection with the cosets of $N_{red}(D^*)$ in $k^*$. In case $(ii)$ the
vertices lying above $v_*$ are in bijection with the (multiplicative!) cosets
of the ``reduced norm'' $N_{red}(k(\sqrt{\rho},\sqrt{\sigma})^*)$ in $k^*$.
(Note that every element of $k(\sqrt{\rho},\sqrt{\sigma})$ is contained in
a quadratic extension of $k$.)
\\ \\
{\bf Remark 3.1.}
Since $N_{k(\overline{t})/k}(k(\overline{t})^*)$ contains all squares
from $k^*$, we see that the group
$k^*/N_{k(\overline{t})/k}(k(\overline{t})^*)$
has exponent $2$ (or $1$). Hence the same holds for its subgroups and
factor groups.
This means that the number of vertices of $\Gamma\backslash\TTT$ lying
above $v_*$ is a power of $2$ (including $1=2^0$ and $\infty$), and
likewise for the number of edges attached to each such vertex.
\\ \\
\noindent {\bf Lemma 3.2.} \it \begin{itemize}
\item[(i)] Let $(G\backslash\TTT)^*$ be the contraction of
$G\backslash\TTT$ obtained by deleting $v_*$ and $\epsilon_*$. Then
$$\phi^{-1}((G\backslash\TTT)^*) \cong (G\backslash
\TTT)^*.$$
\item[(ii)] $$\Gamma \backslash \TTT\; is\; a\; tree
\Longleftrightarrow \det(S_{v_*})=\det(S_{\epsilon_*}).$$
\item[(iii)] $$ \Gamma\backslash\TTT \cong
G\backslash\TTT\Longleftrightarrow
\det(S_{v_*})=\det(S_{\epsilon_*})=k^*.$$
\item[(iv)] Let $\mu=|\det(S_{v_*}):\det(S_{\epsilon_*})|$. If $\mu$ is finite there exists an
epimorphism
$$ \theta: \Gamma\twoheadrightarrow F_{\mu-1},$$
where $F_{\mu-1}$ is the free group of rank $\mu-1$. \noindent If
$\mu$ is infinite of cardinality $\cc$, then there exists an
epimorphism
$$ \theta: \Gamma\twoheadrightarrow F_{\cc},$$
where $F_{\cc}$ is the free group of rank $\cc$.
\end{itemize}
\noindent {\bf Proof.} \rm Let $ w \in \{
v_n,\epsilon_n: n \geq 0\}$. Then $|\phi^{-1}(w)|=1$ by Lemma 2.3.
Part (i) follows. Parts (ii), (iii) are clear. \noindent For part
(iv) the theory of groups acting on trees [Se2, p.42] shows that
$\Gamma$ maps onto the fundamental group $\pi_1(\Gamma\backslash
\TTT$) of the graph $\Gamma\backslash\TTT$. \hfill $\Box$
\\ \\
\noindent {\bf Example 3.3: Hamilton's quaternions.}  Consider the
case where $k=\RR$ and $\rho=\sigma=1$ (so that $\CCC=\RR[x,y]$,
where $x^2+y^2=-1$). Here $D$ consists of {\it Hamilton's
quaternions}. From the proofs of Lemmas 2.5 and 2.7 it follows that
$\det(S_{v_*})=\det(S_{\epsilon_*})=\RR^{+}=(0,\infty)$. Hence by
Lemma 3.2 the fundamental domain $\Gamma\backslash\TTT$ is
isomorphic to the following subtree of $\TTT$

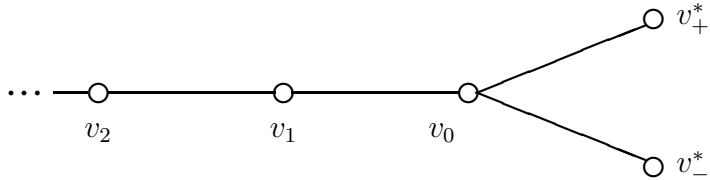
\begin{center}
\setlength{\unitlength}{1pt} \thicklines
\begin{picture}(300,50)(-150,-30)
\put(-110,0){\circle{7}} \put(-40,0){\circle{7}}
 \put(30,0){\circle{7}} \put(100,28){\circle{7}}
\put(100,-28){\circle{7}}\put(-106,0){\line(1,0){63}}
\put(-37,0){\line(1,0){63}}
 \put(-127,0){\line(1,0){13}}
 \put(33,0){\line(5,2){64}}
 \put(33,0){\line(5,-2){64}}\put(-133,0){\circle{.5}}
\put(-138,0){\circle{.5}} \put(-143,0){\circle{.5}}
\put(-110,-15){\makebox(0,0){$v_{2}$}}
\put(-40,-15){\makebox(0,0){$v_1$}}
\put(20,-15){\makebox(0,0){$v_0$}}
\put(115,28){\makebox(0,0){$v^*_{+}$}}
\put(115,-28){\makebox(0,0){$v^*_{-}$}}
\end{picture}
\end{center}

\noindent The vertex $v_*=v(\pi,t)$ lifts (under $\phi$) to the
vertices $v_*^{+}=v(\pi,t)$ and $v_*^{-}=v(\pi,-t)$ (where
$t=y/x$). The edge $\epsilon_*$ lifts in a similar way.
\\ \\
\noindent {\bf Example 3.4: Local constant fields.}
Consider the case where $k$ is a {\it local field}. Thus $k$ is the
field $\QQ_p$ of $p$-adic numbers or a finite extension of $\QQ_p$
or a Laurent series field $\FF_q((u))$.
\par
First we treat the case where the characteristic is different from $2$.
From the proofs of Lemmas 2.5 and 2.8
$$\det(S_{v_*})=\{\alpha^2+\beta^2\rho+\sigma(\gamma^2+\delta^2\rho):(\alpha,\beta,\gamma,\delta)\neq(0,0,0,0)\}$$
and
$$\det(S_{\epsilon_*})=\{\alpha^2+\beta^2\rho:(\alpha,\beta) \neq
(0,0)\}.$$ \noindent For each $\tau \in k$, for which $\sqrt{-\tau}
\notin k$, let
$$N(\tau)=\{\alpha^2+\beta^2\tau: (\alpha,\beta)\neq (0,0)\}.$$
\noindent Then from standard results in the theory of local fields
[Se1, Corollary, p.172]
$$|k^*:N(\tau)|=2.$$
\noindent It follows that
$$|\det(S_{v_*}):\det(S_{\epsilon_*})| \leq 2.$$
\noindent Suppose now that $\det(S_{v_*})=\det(S_{\epsilon_*})$.
Then $N(\rho)=N(\sigma)$ and so $k(\sqrt{-\rho})=k(\sqrt{-\sigma})$,
by [Se1, Proposition 4, p.172]. Thus $\rho=\lambda^2\sigma$, for
some $\lambda \in k$. (Recall that the characteristic of $k$ is not
$2$.) We may therefore assume that $\rho=\sigma$, in which case
$\sqrt{-1} \notin k$. \noindent Now every element in a finite field
is the sum of $2$ squares. Since $\mathrm{char}(k) \neq 2$ we deduce
from this that there exist $\zeta,\eta \in k$ such that
$$\zeta^2+\eta^2+\rho=0.$$
\noindent Then, using elements in $k(\sqrt{-1})$ we can find
$\lambda,\mu \in k$ such that
$$\lambda^2+\rho(1+\mu^2)=0,$$
\noindent which contradicts the properties of $\rho$. \noindent
Hence $\det(S_{v_*})=k^*$ and
$|\det(S_{v_*}):\det(S_{\epsilon_*})|=2$. By Lemma 3.2 therefore
$\Gamma_0\backslash\TTT$ is then isomorphic to the following.

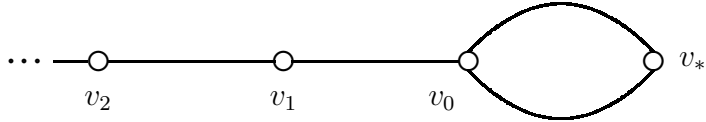
\begin{center}
\setlength{\unitlength}{1pt} \thicklines
\begin{picture}(300,50)(-150,-30)
\put(-110,0){\circle{7}} \put(-40,0){\circle{7}}
 \put(30,0){\circle{7}} \put(100,0){\circle{7}}
\put(-106,0){\line(1,0){63}} \put(-37,0){\line(1,0){63}}
 \put(-127,0){\line(1,0){13}} \put(-133,0){\circle{.5}}
\put(-138,0){\circle{.5}} \put(-143,0){\circle{.5}}
\put(-110,-15){\makebox(0,0){$v_{2}$}}
\put(-40,-15){\makebox(0,0){$v_1$}}
\put(20,-15){\makebox(0,0){$v_0$}} \put(115,0){\makebox(0,0){$v_*$}}
\qbezier(30,3.5)(65,40)(100,3.5)\qbezier(30,-3.5)(65,-40)(100,-3.5)
\end{picture}
\end{center}
\noindent As in Example 3.3 the edge $\epsilon_*$ lifts to a pair of
edges in $\Gamma\backslash \TTT$. However the vertex $v_*$ lifts to
only one. By Lemma 3.2 (iv) there exists an epimorphism
$$\theta: \Gamma\twoheadrightarrow \ZZ.$$

\noindent Now assume that the characteristic is $2$. Then by Example 1.5
we are in case $(iii)$ or $(iv)$ from Theorem 1.1. As $N_{red}(D^*)=k^*$
by [R, Theorem 33.4], there is one vertex of $\Gamma\backslash \TTT$ above
$v_*$. In case $(iv)$ the norm of $k(\overline{t})^*$ has index $2$ in
$k^*$ and the graph is the same as in the other characteristics. But in
case $(iii)$ the norm from $k(\sqrt{\rho})^*$ is surjective and hence the
graph is the same as for $G$ (see Theorem 2.9).
\\ \\
\noindent {\bf Example 3.5: Global constant fields.}
Now we investigate the case where $k$ is a global field, i.e. an
algebraic number field or an algebraic function field of one
variable over a finite constant field.
This will also show the extent to which the structure of
$\Gamma\backslash \TTT$ varies with the field $k$.
\\ \\
{\bf Proposition 3.6.} \it
If $k$ is an algebraic number field, the number of vertices of
$\Gamma\backslash\TTT$ above $v_*$ is $2^m$ where $m$ is the number
of real places of $k$ at which the quaternion algebra $D$ does not split.
In other words, $m$ is the number of different embeddings $\varphi$ of
$k$ into $\RR$ for which the equation
$$y^2 +\varphi(\rho)x^2 +\varphi(\sigma)=0$$
has no real solutions. Between each such vertex and $v_0$ there are
infinitely many edges.
\par
If $k$ is a global function field, there is only one vertex above $v_*$.
If moreover $\mathrm{char}(k)=2$ and $\CCC$ is of type $(iii)$ (compare
Theorem 1.1), there is exactly one edge attached to this vertex. In all
other cases there are infinitely many edges attached to it.
\rm
\\ \\
{\bf Proof.} By Example 1.5 the case $(ii)$ cannot occur for a global field.
Thus by Remark 2.7 we have $S(\pi,t)=D^*$ with a non-split quaternion
algebra $D$. By the Theorem of Hasse-Schilling-Maass (see for example
[R, Theorem 33.15]), the reduced norms of $D^*$ in $k^*$ are exactly
the elements that are positive at every real place of $k$ that ramifies
in $D$. Thus $[k^* :N_{red}(D^*)]=2^m$.
\par
If $k$ is a global field of positive characteristic, then $k$ has of course
no real embeddings at all and hence $\Gamma\backslash\TTT$ has only one
vertex above $v_*$.
\par
On the other hand, if $k(\overline{t})$ is separable over $k$, there are
infinitely many places of $k$ that are inert in $k(\overline{t})/k$, and
therefore in $k^*$ the norms of $k(\overline{t})$ are subject to infinitely
many congruence conditions. So $N_{k(\overline{t})/k}(k(\overline{t})^*)$
has infinite index in $k^*$ and hence also in $N_{red}(D^*)$.
\par
If $k(\overline{t})/k$ is inseparable however, the norm is surjective.
\hfill $\Box$
\\ \\
\noindent Proposition 3.6 could also have been proved using the
theory of quadratic forms.
\\ \\
\noindent {\bf Example 3.7: Constant fields with characteristic 2.}
\noindent Now we show that for case $(ii)$ the structure of
$\Gamma\backslash \TTT$ depends only on the constant field $k$.
\\ \\
\noindent {\bf Lemma 3.8.} \it Let $\mathrm{char}(k)=2$ and suppose
that $\CCC$ is of type $(ii)$ (compare Theorem 1.1). Then:
\begin{itemize}
\item[a)] If $[k:k^2]=4$, then there is exactly one vertex in
$\Gamma\backslash \TTT$ above $v_*$.
\item[b)] If $[k:k^2]>4$, then there are infinitely many vertices
above $v_*$.
\item[c)] In either case there are infinitely many edges between each
such vertex and $v_0$.
\end{itemize}
\rm
\noindent {\bf Proof.}\\
a) If $[k:k^2]=4$, then $k^2 +\rho k^2 +\sigma k^2 +\rho\sigma
k^2=k$, so $N_{red}(k(\sqrt{\rho},\sqrt{\sigma})^*)=k^*$.
\par
b) In this case $k$ is a vector space (of dimension $>1$) over the
field $E=k^2 +\rho k^2 +\sigma k^2 +\rho\sigma k^2$. So $k^*/E^*$ is
a projective space over $E$. Since $E$ is infinite, this projective
space has infinitely many points.
\par
c) For the same reason as in b) the index of
$N_{k(\sqrt{\rho})/k}(k(\sqrt{\rho})^*)$ in
$N_{red}(k(\sqrt{\rho},\sqrt{\sigma})^*)$ is infinite. (The quotient
forms a projective line over the field $k^2 +\rho k^2$.)
\hfill $\Box$
\\ \\
\noindent {\bf Example 3.9: Laurent series fields as constant fields.}
If the characteristic of the field $L$ is different from $2$, then
a nonzero element from $k=L((u))$ is a square in $k$ if and only if its
first nonzero coefficient is a square in $L^*$ and the corresponding
exponent of $u$ is even. In particular, if $L$ is closed under taking
square-roots, then $k^*/(k^*)^2$ is represented by $1$ and $u$, and hence
one cannot find $\rho, \sigma\in k$ as in Theorem 1.1. So if $L$ is closed
under taking square-roots, there cannot be a nonrational genus zero
function field with constant field $L((u))$.
\par
We now restrict our attention to Laurent series fields over $\RR$
and leave it as an easy exercise for the reader to verify the following
examples, all over $k=\RR((u))$:
$$\begin{array}{ll}
y^2 + x^2 +1=0,\ & \hbox{\rm $4$ vertices above $v_*$, each one with $1$ edge;}\\
y^2 + x^2 +u=0,\ & \hbox{\rm $2$ vertices above $v_*$, each one with $2$ edges;}\\
y^2 +ux^2 +1=0,\ & \hbox{\rm $2$ vertices above $v_*$, each one with $1$ edge.}\\
\end{array}$$
This shows that over the same constant field the structure of
$\Gamma\backslash\TTT$ may vary with $\CCC$. Note that in the last two
instances even the fields $K$ are isomorphic and the associated quaternion
algebras and quaternion groups are also isomorphic.
\\ \\
\noindent To conclude we record one of the most elementary cases for which
$\Gamma\backslash \TTT$ has infinitely many edges.
Let $k =\QQ$ and $\rho=\sigma=1$. In this case
$\CCC=\mathcal{Q}=\QQ[x,y]$, where $x^2+y^2=-1$. \noindent From the
{\it Four-Square Theorem} it follows that
$$S_{v_*}=\QQ^{+}.$$
\noindent On the other hand let $p_i$ be a rational prime where $p_i
\equiv 3\;(\mod 4)\;(1\leq i\leq t)$ Then
$$p_1\cdots p_t \notin S_{\epsilon_*}.$$
\noindent In this case therefore the index
$|S_{v_*}:S_{\epsilon_*}|$ is infinite. By Lemma 3.2 the graph of
$\Gamma\backslash \TTT$ is then of the following form.
\\ \\

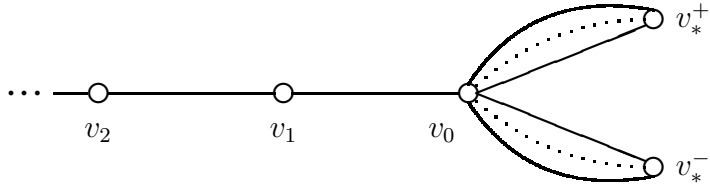
\begin{center}
\setlength{\unitlength}{1pt} \thicklines
\begin{picture}(300,50)(-150,-30)
\put(-110,0){\circle{7}} \put(-40,0){\circle{7}}
 \put(30,0){\circle{7}} \put(100,28){\circle{7}}
\put(100,-28){\circle{7}}\put(-106,0){\line(1,0){63}}
\put(-37,0){\line(1,0){63}}
 \put(-127,0){\line(1,0){13}}
 \put(33,0){\line(5,2){64}}
 \put(33,0){\line(5,-2){64}}\put(-133,0){\circle{.5}}
\put(-138,0){\circle{.5}} \put(-143,0){\circle{.5}}
\put(-110,-15){\makebox(0,0){$v_{2}$}}
\put(-40,-15){\makebox(0,0){$v_1$}}
\put(20,-15){\makebox(0,0){$v_0$}}
\put(115,28){\makebox(0,0){$v_*^{+}$}}
\qbezier(30,3.5)(54,40)(100,31.5) \qbezier[15](33,3)(60,27)(96.5,28)
\qbezier(30,-3.5)(54,-40)(100,-31.5)
\qbezier[15](33,-3)(60,-27)(96.5,-28)\put(115,-28){\makebox(0,0){$v_*^{-}$}}
\end{picture}
\end{center}
\noindent The dotted lines mean that there is a countably infinite number
of edges between the corresponding vertices.
\\ \\
\noindent From Lemma 3.2 (iv) we conclude with the following.
\\ \\
\noindent {\bf Corollary 3.10.} \it With the above notation, there
exists an epimorphism
$$\theta: SL_2(\mathcal{Q}) \twoheadrightarrow F_{\omega},$$
where $F_{\omega}$ is the free group of countably infinite rank.

\noindent \rm Every finitely (or countably) generated group
(including every finite group) is therefore an image of
$SL_2(\mathcal{Q})$.
\\ \\

\subsection*{\hspace*{10.5em} References}
\begin{itemize}
\item[{[A]}] E. Artin, \it Algebraic Numbers and Algebraic Functions,
\rm Gordon and Breach, New York, 1968.
\item[{[FJ]}] M. Fried and M. Jarden, \it Field Arithmetic (Second Edition),
\rm Springer, Berlin, Heidelberg, New York, 2004.
\item[{[G]}] E.-U. Gekeler, \it Drinfeld Modular Curves, \rm Lecture Notes
in Mathematics, vol. 1231, Springer, Berlin, Heidelberg, New York, 1986.
\item[{[GMV]}] F. Grunewald, J. Mennicke and L. Vaserstein, On the
groups $SL_2(\ZZ[x])$ and $SL_2(k[x,y])$, {\it Israel J. Math}
\bf{86} \rm(1994), 157-193.
\item[{[KM]}] S. Krstic and J. McCool, Free quotients of
$SL_2(R[x])$, {\it Proc. Amer. Math. Soc.} \bf{125} \rm(1997),
1585-1588.
\item[{[L]}]  A. Lubotzky, Lattices in rank one Lie groups over local fields,
Geom. Funct. Anal. \bf{1} \rm (1991), 405-431.
\item[{[M1]}] A.W. Mason, Serre's generalization of Nagao's theorem: an elementary
approach, {\it Trans. Amer. Math. Soc.} \bf{353} \rm (2003), 749-767.
\item[{[M2]}] A. W. Mason, The generalization of Nagao's
theorem to other subrings of the rational function field, \it
Comm. Algebra \bf{31} \rm(2003), 5199-5242.
\item[{[MS]}] A. W. Mason and A. Schweizer, The minimum index
of a non-congruence subgroup of $SL_2$ over an arithmetic domain,
\it Israel J. Math. \bf 133 \rm (2003), 29-44.
\item[{[R]}] I. Reiner, \it Maximal orders, \rm Academic Press,
London - New York, 1975.
\item[{[Se1]}] J.-P. Serre, \it Local Fields, \rm Springer, Berlin,
Heidelberg, New York, 1979.
\item[{[Se2]}] J.-P. Serre, \it Trees, \rm Springer, Berlin,
Heidelberg, New York, 1980.
\item[{[St]}] H. Stichtenoth, \it Algebraic Function Fields and Codes,
\rm Springer, Berlin, Heidelberg, New York, 1993.
\item[{[T]}] S. Takahashi, The fundamental domain of the
tree of $GL(2)$ over the function field of an elliptic curve, \it
Duke Math. J. \bf 72 \rm (1993), 85-97.
\item[{[VS]}] G. Villa Salvador, \it Topics in the Theory of Algebraic
Function Fields, \rm Birkh\"auser, Boston, 2006
\end{itemize}
\end{document}